\documentclass{amsart}

\usepackage{amscd,amssymb,amsthm,amsmath,enumerate,textcomp} 

\usepackage[utf8]{inputenc}
\usepackage[T1]{fontenc}
\usepackage[svgnames,hyperref]{xcolor}
\usepackage{lmodern}
\usepackage{xspace} 
\usepackage[nomessages]{fp} 
\usepackage{ifthen} 
\usepackage[normalem]{ulem}
\usepackage{comment}

\usepackage{tikz}
\usetikzlibrary{cd,arrows}
\tikzset{>=stealth}
\tikzcdset{arrow style=tikz}
\tikzset{link/.style={column sep=1.8cm,row sep=0.16cm}}
\usepackage{xy}

\newcommand{\mapsrat}{ \, \tikz[baseline=-.6ex] \draw[|->,densely dashed,line width=.5] (0,0) -- +(.5,0); \, }
\newcommand{\rat}{ \, \tikz[baseline=-.6ex] \draw[->,densely dashed,line width=.5] (0,0) -- +(.5,0); \, }
\renewcommand{\to}{ \, \tikz[baseline=-.6ex] \draw[->,line width=.5] (0,0) -- +(.5,0); \, }

\usepackage[shortlabels]{enumitem}
\setlist{wide}
\setlist[enumerate]{label=\rm{(\arabic*)}}
\setlist[enumerate,2]{label=\rm({\it\roman*})}
\setlist[itemize]{label=\raisebox{0.25ex}{\tiny$\bullet$}}

\newcommand{\nicecolor}{Navy}
\usepackage[pdfauthor={H. Bergner \& S. Zimmermann}, backref=page, colorlinks, linktocpage, citecolor=\nicecolor, linkcolor=\nicecolor, urlcolor=\nicecolor]{hyperref}
\usepackage[all]{hypcap} 

\newcommand{\p}{\mathbb{P}}
\newcommand{\A}{\mathbb{A}}
\newcommand{\Z}{\mathbb{Z}}
\newcommand{\N}{\mathbb{N}}
\newcommand{\C}{\mathbb{C}}
\newcommand{\R}{\mathbb{R}}
\newcommand{\Q}{\mathbf{Q}}
\renewcommand{\k}{\mathbf{k}}
\DeclareMathOperator{\Bir}{Bir}
\DeclareMathOperator{\Aut}{Aut}
\newcommand{\id}{\text{\rm id}}
\newcommand{\ord}{\mathrm{ord}}
\newcommand{\bR}{\mathrm{R}}

\theoremstyle{plain}
\newtheorem{theorem}{Theorem}[section]
\newtheorem{corollary}[theorem]{Corollary}
\newtheorem{proposition}[theorem]{Proposition}
\newtheorem{lemma}[theorem]{Lemma}

\theoremstyle{definition}
\newtheorem{definition}[theorem]{Definition}
\newtheorem{question}[theorem]{Question}
\newtheorem{example}[theorem]{Example}
\newtheorem{remark}[theorem]{Remark}

\title[Properties of Cremona groups in the Euclidean topology]{Properties of the Cremona group endowed with the Euclidean topology}

\author{Hannah Bergner and Susanna Zimmermann}

\address{Hannah Bergner}
\email{hannah.bergner-c9q@rub.de}

\address{Susanna Zimmermann, Univ Angers, CNRS, LAREMA, SFR MATHSTIC, F-49000 Angers, France}
\email{susanna.zimmermann@univ-angers.fr}


\thanks{
During this project, the second author was supported by the ANR Project FIBALGA ANR-18-CE40-0003-01, the Projet PEPS 2019 ``JC/JC'', the Project \'Etoiles montantes of the R\'egion Pays de la Loire and the Centre Henri Lebesgue, program ANR-11-LABX-0020-0. 
}

\begin{document}

\begin{abstract}
Consider a Cremona group endowed with the Euclidean topology introduced by Blanc and Furter. It makes it a Hausdorff topological group that is not locally compact nor metrisable.
We show that any sequence of elements of the Cremona group of bounded order that converges to the identity is constant. 
We use this result to show that Cremona groups do not contain any non-stationary sequence of subgroups converging to the identity. We also show that, in general, paths in a Cremona group do not lift and do not satisfy a property similar to the definition of morphisms to a Cremona group. 
\end{abstract}

\maketitle

\section{Introduction}

The Cremona group $\Bir(\p^n_{\k})$ denotes the group of birational transformations of~$\p^n_{\k}$ defined over a field $\k$. 
If $\k$ is a local field, that is, a locally compact topological field with respect to a non-discrete topology, it can be endowed with the so-called Euclidean topology (for a precise definition see Section~\ref{sec:topologies}), which comes from the Euclidean topology of the coefficients of birational maps, and which gives $\Bir(\p^n_{\k})$ the structure of a Hausdorff topological group which is not locally compact and not metrisable. Moreover, the restriction of the Euclidean topology to any algebraic subgroup of $\Bir(\p^n_{\k})$ is the usual Euclidean topology on that algebraic group \cite[Theorem 3, Lemma 5.15, Lemma 5.16]{BF13}.
While $\Bir(\p^2_\k)$ is compactly presented in the Euclidean topology by a quadratic involution of $\p^2_\k$ and a compact subset of $\Aut(\p^2)$ \cite[Theorem A]{Z16}, the group $\Bir(\p^n_\k)$, $n\geq3$ is not generated by a compact subset (Proposition~\ref{lem:properties_topology}\ref{pro_top:7}). In fact, $\Bir(\p^3_\C)$ is not even generated by its algebraic subgroups \cite[Theorem C]{BY19}. 
The Euclidean topology on Cremona groups is largely unstudied and results can be found in \cite{BF13,BZ16,UZ21,Z16}. 

This article is motivated by the question whether $\Bir(\p^n_\k)$ has no small subgroups. A topological group $G$ is said to have no small subgroups, if there exists a neighbourhood of the identity element that contains no non-trivial subgroups of $G$. For instance, any complex or real Lie group has no small subgroups, see for instance \cite[Proposition 2.17]{HM07}. 
The property of having/not having small subgroups is related to Hilbert's $5$th problem: 
any finite dimensional, locally compact, separable metric, locally connected group without small subgroup is a (complex or real) Lie group \cite[Theorem 3.1]{G52}. 
It implies that any finite dimensional, locally compact, separable metric, locally connected groups is a generalised (complex or real) Lie group \cite[p.440]{MZ52}. Both results were generalised to the infinite dimensional case in \cite[Theorem 5]{Yamabe}.

If $\Bir(\p^n_{\k})$ is endowed with the Euclidean topology, it is not metrisable and it follows from \cite[Theorem 2.1.1]{E70} that $\Bir(\p^n_{\k})$ is not uniformly free from small subgroups. 
However, this does not mean that $\Bir(\p^n_{\k})$ has small subgroups. 
If one wants to find small subgroups, a first intuitive idea could be to construct a sequence of involutions that converge to the identity. 
However, a first example of a group containing small subgroups is $\mathrm{SL}_2(\Z_p)$, whose small subgroups contain no finite order elements. Also for the Cremona group, the intuitive approach fails:

\begin{theorem} \label{thm1}
Let $\k$ be a local field of characteristic zero. 
Let $n,D\geq2$ and $ f_m\in\Bir(\p^n_\k)$ be a sequence with $\ord(f_m)\le D$ for $m\gg0$ converging to $\id\in\Bir(\p^n_\k)$ in the Euclidean topology. 
Then we have $f_m=\id$ for $m\gg 0$.
\end{theorem}

This is similar to a general property of analytic groups over an ultrametric field~$\k$: if $G$ is such a group, then there is an open subgroup $U\subset G$ that contains no finite subgroup whose order is prime to $\mathrm{char}(\k)$ \cite[Chapter IV, Remark after Theorem~5, LG~4.27]{Serre-LA+LG}. 

Let $(G_m)_{m\in \N}$ be a sequence  of subgroups of a topological group $G$. 
 We say that it converges to the neutral element $\id\in G$ if for any open neighbourhood~$U$ of $\id$ there exists $m_0\in \N$ such that $G_m \subset U$ for all $m> m_0$.  If $G$ is a topological group containing a sequence of subgroups $(G_m)_{m\in\N}$ that converges to $\id$ and such that $G_m\neq\{\id\}$ for each $m\geq1$, then $G$ has small subgroups and for first-countable topological groups the converse also holds true.
Note however that $\Bir(\p^n_\k)$ is not first-countable (cf. \cite[Lemma 5.16 and its proof]{BF13}) and  might have small subgroups even if there does not exist any non-trivial sequence of subgroups $(G_m)_{m\in\N}$ converging to the identity.
 
Using Theorem~\ref{thm1}, we show the following result.

\begin{theorem}\label{theorem:converging subgroups}
Let $\k=\R$ or $\k=\C$. 
  Let $(G_m)_{m\in \N}$ be a sequence of subgroups of $\Bir(\p^n_{\k})$ that converges to the identity $\id \in \Bir(\p^n_{\k})$.
  Then we have $G_m=\{\id\}$ for $m\gg 0$.
\end{theorem}

Theorem~\ref{theorem:converging subgroups} fails over the p-adic numbers, see Remark~\ref{rmk:fail for p-adic}. 

The following is a corollary of Theorem~\ref{thm1}. See Definition~\ref{Defi:Morphisms} for the definition of a morphism from a variety to $\Bir(\p^n)$. 

\begin{corollary}\label{thm2}
Let $\k$ be a local field of characteristic zero. 
Let $n,D\geq2$ and $\rho\colon \A^1_{\k}\to \Bir(\p^n_\k)$ be a morphism from the variety $\A_\k^1$ to $\Bir(\p^n_\k)$ with $\rho(0)=\id$. Suppose that there is a sequence $(t_m)_m\subset\A^1_{\k}$ converging to $0$ in the Euclidean topology with $t_m\neq 0$ for all $m$ and such that $\ord(\rho(t_m))\leq D$ for $m\gg0$. 
Then $\rho(\A^1_{\k})=\{\id\}$.
\end{corollary}

Along the way of proving the main results, we prove in Corollary~\ref{cor:compact_subgroups} that for $\k=\C$ or $\k=\R$ every compact subgroup of $\Bir(\p^n_\k)$ carries a compatible {\em real} Lie group structure. \smallskip

An open problem is to describe the fundamental group of the topological group $\Bir(\p^n_\k)$. 
Let us recall the definition of the Euclidean topology on $\Bir(\p^n_\k)$. 
For any $d\geq1$, let $\Bir(\p^n_\k)_{\leq d}\subset\Bir(\p^n_\k)$ be the subset of elements of degree $\leq d$. 
The Euclidean topology on $\Bir(\p^n_\k)$ is defined as the inductive topology of the sets $\Bir(\p^n_\k)_{\leq d}$, $d\geq1$, each endowed with the Euclidean topology, which is defined as follows: consider the vector space $(\k[x_0,\dots,x_n]_d)^{n+1}$ of $n+1$-tuples of homogeneous polynomials of degree $d$ in $n+1$ variables. A non-zero tuple $(f_0,\dots,f_n)$ induces the rational map $\p^n_\k\rat\p^n_\k$, $[x_0:\cdots:x_n]\mapsrat[f_0(x_0,\dots,x_n):\dots:f_n(x_0,\dots,x_n)]$. Consider the subset $H_d(\k)\subset\p(\k[x_0,\dots,x_n]_d)$ of elements that induce a birational map of $\p^n_\k$. There is a natural surjective map $\pi_d\colon H_d(\k)\to\Bir(\p^n_\k)_{\leq d}$. The Euclidean topology on $\Bir(\p^n_\k)_{\leq d}$ is defined to be the quotient topology with respect to $\pi_d$. 

If $\Omega$ is a compact set, then the image of a continuous map $\rho\colon \Omega\to\Bir(\p^n_\k)$ is contained in some $\Bir(\p^n_\k)_{\leq d}$, see Proposition~\ref{lem:properties_topology}. 
A first idea to study the fundamental group of $\Bir(\p^n_\k)$ could be to lift a loop $\rho$ to a loop in $H_d(\k)$. We show
however that this fails in general.

\begin{proposition}
\label{prop:nolifting-main}
Let $\k=\R$ or $\k=\C$. 
For any $n\geq2$, there exists a map $\rho\colon [-1,1]\to\Bir(\p^n_\k)$ that is continuous in the Euclidean topologies such that for 
arbitrarily small $\varepsilon>0$ and any $d\geq1$, there is no continuous map $\rho_{\varepsilon}\colon(-\epsilon,\epsilon)\to H_d(\k)$ with $\rho|_{(-\epsilon,\epsilon)}=\pi_d\circ\rho_{\varepsilon}$.
\end{proposition}

The key ingredient for the construction of the examples is the existence of space filling curves of $\mathrm{PSU}(n+1)$. Similar examples can be constructed over non-Archimedean fields, see Remark~\ref{rmk:Cantor}. 

The Euclidean topology is a refinement of the so-called Zariski topology. The Zariski topology can be defined via the concept of morphisms from algebraic varieties to the Cremona group, see Section~\ref{sec:topologies}. We show that Euclidean-continuous maps $[0,1]\to\Bir(\p^n_\k)$ do not satisfy in general a condition similar to the definition of morphisms of algebraic varieties to $\Bir(\p^n_\k)$:

\begin{proposition}
\label{prop:example1-main}
Let $\k=\R$ or $\k=\C$. 
For every $n\geq2$ there exists a map $\rho\colon [0,1]\to\Bir(\p^n_\k)$ that is continuous in the Euclidean topologies such that there is no open non-empty subset $U\subset[0,1]\times\p^n_\k$ such that $U\to[0,1]\times\p^n_\k$, $(t,p)\mapsto(t,\rho(t)(p))$ is a well defined map.
\end{proposition}

\bigskip

\noindent{\bf Acknowledgements:} The authors thank Jérémy Blanc and Christian Urech for fruitful discussions early on in the project, J\"org Winkelmann for helpful discussion about the proof of Proposition~\ref{prop:sequences_of_bounded_order}, the referee for suggesting an approach over the $p$-adic numbers and Philipp Habegger for his helpful remarks on non-Archimedean fields.

\section{Topologies on the Cremona groups}
\label{sec:topologies}
Let $\k$ be any field. 
The Zariski topology is defined over any field and the Euclidean topology over any local field. 

\subsection{The Zariski topology}\label{ss:zariskitopology}

The notion of morphism from an algebraic variety $A$ to $\Bir(X)$ was introduced by M.~Demazure in \cite{D70}. 

\begin{definition}\label{Defi:Morphisms}
Let $X$ be an irreducible algebraic variety and $A$ an algebraic variety. Consider a birational map of the form $A\times X\dashrightarrow A\times X$, $(a,x)\dashrightarrow(a,f_a(x))$, that induces an isomorphism $U\to V$ on open dense subsets $U,V\subset A\times X$ whose projection onto $A$ is surjective. 
The map $A(\k)\ni a\mapsto f_a$ represents a map from $A(\k)$ to $\Bir(X)$, and will be called a \emph{morphism} from $A$ to $\Bir(X)$. It is denoted by $A\to \Bir(X)$.
\end{definition}

Even if $\Bir(X)$ is not representable by an algebraic variety or an ind-algebraic variety if $X$ is rational \cite{BF13}, 
we can define a topology on the group $\Bir(X)$ compatible with morphisms. This topology is called \emph{Zariski topology} by J.-P. Serre in \cite{Se}:

\begin{definition}  \label{defi: Zariski topologyk}
Let $X$ be an irreducible algebraic variety. A subset $F\subseteq \Bir(X)$ is defined to be \emph{closed in the Zariski topology} if for any algebraic variety $A$ and any morphism $A\to \Bir(X)$ the preimage of $F$ in $A$ is closed.
\end{definition}

In this definition one can of course replace ``any algebraic variety $A$'' with ``any \emph{irreducible} algebraic variety $A$''.

\begin{definition}\label{DefWHG}\cite[Definition 2.3]{BF13}
Let $d,n$ be positive integers. 
\begin{enumerate}
\item

We define $W_d(\k)$ to be the projective space parametrising equivalence classes of non-zero  $(n+1)$-tuples $(h_0,\dots,h_n)$
of homogeneous polynomials $h_i\in \k[x_0,\dots,x_n]$ of degree $d$,
where $(h_0,\dots,h_n)$ is equivalent to $(\lambda h_0,\dots,\lambda h_n)$ for any $\lambda\in \k^{*}$.
The equivalence class of $(h_0,\dots,h_n)$ will be denoted by $[h_0:\dots:h_n]$.
\item
We define $H_d(\k)\subseteq W_d(\k)$ to be the set of elements $h=[h_0:\dots:h_n]\in W_d(\k)$
such that the rational map
$\psi_h\colon \p^n_{\k}\dasharrow \p^n_{\k}$ given by 
\[[x_0:\dots:x_n]\dasharrow
[h_0(x_0,\dots,x_n):\dots:h_n(x_0,\dots,x_n)]\] is birational.
We denote by $\pi_d$ the map $H_d(\k)\to \Bir(\p^n_{\k})$ which sends $h$ onto~$\psi_h$.
\item We denote by $W_d$ the underlying projective space over $\k$ whose set of rational points is $W_d(\k)$ and by $H_d$ the subset whose set of rational points is $H_d(\k)$.
\end{enumerate}
\end{definition}

\begin{proposition}\label{Prop:HdPidKMorphism}
Let $d,n$ be positive integers and $\k$ any field.
\begin{enumerate}
\item
The set $H_d$ is locally closed in the projective space $W_d$ and thus inherits the structure of an algebraic variety;
\item
The map $\pi_d$ corresponds to a morphism $H_d\to \Bir(\p^n_{\k})$ and its image is the set $\Bir(\p^n_{\k})_{\le d}$ of all birational maps of degree $\le d$.
\item 
If $\k\subset L$ is a field extension, there are canonical inclusions 
$\hat\iota\colon H_d(\k)\hookrightarrow H_d(L)$ and 
$\iota\colon\Bir(\p^n_\k)_{\leq d}\hookrightarrow\Bir(\p^n_L)_{\leq d}$ and associated maps $(\pi_d)_{\k}\colon H_d(\k)\to\Bir(\p^n_\k)_{\leq d}$ and $(\pi_d)_L\colon H_d(L)\to\Bir(\p^n_L)_{\leq d}$.  
Then 
\[
\iota\circ(\pi_d)_{\k}=(\pi_d)_L\circ\hat\iota
\quad\text{and}\quad 
(\pi_d)_{\k}^{-1}(A)=(\pi_d)_{L}^{-1}(A)\cap H_d(\k)
\]
for any set $A\subset\Bir(\p^n_\k)$. 
\end{enumerate}
\end{proposition}
\begin{proof}
The first two claims follows from \cite[Lemma 2.4]{BF13}. 
The first part of the third claim is straight forward to check. For the last claim, it suffices to remark the following. 
If $h\colon[x_0:\dots:x_n]\mapsto[h_0(x_0,\dots,x_n):\dots:h_n(x_0,\dots,x_n)]$ is contained in $\Bir(\p^n_L)_{\leq d}$, then its pre-image $(\pi_d)_L^{-1}(h)$ in $H_d(L)$ is the set of elements $[Hh_0:\dots:Hh_n]\in H_d(L)$ such that $H\in L[x_0,\dots,x_n]$ is homogeneous of degree $d-\deg(h)$. 
The set $(\pi_d)_{\k}^{-1}(h)$ in $H_d(\k)$ is the set of elements $[Hh_0:\dots:Hh_n]\in H_d(\k)$ where $H\in \k[x_0,\dots,x_n]$ is homogeneous of degree $d-\deg(h)$. 
\end{proof}

\subsection{The Euclidean topology}
Let $\k$ be a local field, that is, a locally compact topological field with respect to a non-discrete topology. Any local field of characteristic zero is isomorphic to $\R$, $\C$ or a finite extension of the $p$-adic number $\Q_p$. Any local field of positive characteristic $p>0$ is isomorphic to the field $\mathbb{F}_q(\!(t)\!)$ of Laurent series with $q=p^n$, $n\geq1$. 
The Euclidean topology of $\Bir(\p^n_{\k})$ was first described in \cite[Section 5]{BF13}. It is a refinement of the Zariski topology on $\Bir(\p^n_{\k})$.

\begin{definition}\label{def:Euclidean_top}
The Euclidean topology on $\Bir(\p^n_{\k})$ is defined as follows.
For each $d\ge1$, we endow $W_d(\k)$ 
with the classical Euclidean topology, and $H_d(\k)\subset W_d(\k)$ with the induced Euclidean topology. We define the Euclidean topology on $\Bir(\p^n_{\k})_{\leq d}$ to be the quotient topology given by the map $\pi_d\colon H_d(\k)\to \Bir(\p^n_{\k})_{\le d}$, where $H_d(\k)$ is endowed with the Euclidean topology as just described before. 
The Euclidean topology on $\Bir(\p^n_{\k})$ is defined to be the inductive limit topology induced by the inclusions 
$\Bir(\p^n_{\k})_{\le d}\hookrightarrow \Bir(\p^n_{\k})_{\le d+1}$.
\end{definition}

We will show in Proposition~\ref{prop:counter_example1} that a condition as in Definition~\ref{Defi:Morphisms} of morphisms into Cremona groups does not hold in general for Euclidean continuous maps.

The Euclidean topology on $\Bir(\p^n_{\k})$ has the following properties.

\begin{proposition}\label{lem:properties_topology}
Let $d\geq1$ and $\k$ a local field. We endow $W_d(\k)$ and $H_d(\k)$ with the Euclidean topology and $\Bir(\p^n_{\k})_{\le d}$ with the quotient topology. Then: 
\begin{enumerate}
\item\label{pro_top:1} $\pi_d\colon H_d(\k)\to \Bir(\p^n_{\k})$ is proper and closed.
\item\label{pro_top:2} $\Bir(\p^n_{\k})_{\le d}$ is locally compact and Hausdorff.
\item\label{pro_top:3}  $\Bir(\p^n_{\k})_{\leq d}\hookrightarrow\Bir(\p^n_{\k})_{\le d+1}$ is a closed embedding.  
\item\label{pro_top:4} Any compact subset of $\Bir(\p^n_{\k})$ is contained in some $\Bir(\p^n_{\k})_{\le d}$. In particular, any convergent sequence has bounded degree.
\item\label{pro_top:5}  $\Bir(\p^n_{\k})$ is a Hausdorff topological group.
\item\label{pro_top:6}  For $n\geq2$, $\Bir(\p^n_{\k})$ is not locally compact, not metrisable and not first-countable.
\item\label{pro_top:7} For any $n\geq3$, the group $\Bir(\p^n_{\k})$ endowed with the Euclidean topology is not generated by a compact set.
\end{enumerate}
\end{proposition}
\begin{proof}
These are \cite[Lemmas 5.4, 5.6, 5.8, 5.10, 5.12, 5.13, 5.14, 5.15, 5.16, 5.17]{BF13} in the case of algebraically closed fields of characteristic zero and they can be repeated almost word by word over non-closed fields. 
\end{proof}

If $\k\subset L$ is a field extension, then the Zariski topology on $\Bir(\p^n_\k)$ is finer than the subset topology induced by $\Bir(\p^n_L)$ \cite[Lemma 2.7]{BZ16}.

\begin{lemma}\label{lem:subset-topology}
Let $\k\subset L$ be local fields. 
Then the Euclidean topology on $\Bir(\p^n_\k)$ is finer than the subset topology on $\Bir(\p^n_\k)$ induced by the Euclidean topology on $\Bir(\p^n_L)$. 

In particular, if a sequence $(f_m)_m$ in $\Bir(\p^n_\k)$ converges to $f\in\Bir(\p^n_\k)$ with respect to the Euclidean topology, then it converges to $f$ with respect to the Euclidean topology on $\Bir(\p^n_L)$. 
\end{lemma}
\begin{proof}
If $\hat{A}\subset \Bir(\p^n_\k)$ is closed in the subset topology induced by the Euclidean topology on $\Bir(\p^n_L)$,
then there exists a closed subset $\hat{C}$ in $\Bir(\p^n_L)$ such that $\hat{A}=\hat{C}\cap \Bir(\p^n_\k)$. For any $d\geq1$, the set $A_d=\hat{A}\cap\Bir(\p^n_\k)_{\leq d}$ is closed in $\Bir(\p^n_\k)$ and $C_d=\hat{C}\cap\Bir(\p^n_L)_{\leq d}$ is closed in $\Bir(\p^n_L)_{\leq d}$ and we have $A_d=C_d\cap \Bir(\p^n_\k)_{\leq d}$. 
Proposition~\ref{Prop:HdPidKMorphism} implies that
\[
(\pi_d)_\k^{-1}(A_d)=(\pi_d)_\k^{-1}(C_d\cap \Bir(\p^n_\k)_{\leq d})=(\pi_d)_L^{-1}(C_d)\cap H_d(\k).
\]
The set $(\pi_d)_L^{-1}(C_d)$ is closed in $H_d(L)$ and hence the intersection $(\pi_d)_L^{-1}(C_d)\cap H_d(\k)$ is closed in $H_d(\k)$. Therefore, by Proposition~\ref{lem:properties_topology}\ref{pro_top:1}, the set $A_d$ is closed in $\Bir(\p^n_\k)_{\leq d}$ with respect to the Euclidean topology on $\Bir(\p^n_\k)_{\leq d}$.  
This holds for all $d\geq1$, so $A$ is closed in the Euclidean topology on $\Bir(\p^n_\k)$.
\end{proof}

\section{Sequences of elements of bounded order}
Let $\k$ be a local field. 
Throughout this section, the spaces $\k^n,\p^n_\k,W_d(\k), H_d(\k)$ and $\Bir(\p^n_\k)$ are endowed with the Euclidean topology. We denote by $||\cdot||$ the norm on $\k^n$.

\subsection{Locally uniform convergence}

\begin{remark}\label{rmk:composition}
Let $(f_m)_m\subset\Bir(\p^n_{\k})$ be a sequence of elements converging to $g\in\Bir(\p^n_{\k})$. Since $\Bir(\p^n_{\k})$ is a Hausdorff topological group by Proposition~\ref{lem:properties_topology}\ref{pro_top:5}, it follows that for each $i\geq1$, the sequence of compositions $(f_m^i)_m$ converges to $g^i$. 
\end{remark}

We recall the definition of uniform convergence and locally uniform convergence.

\begin{definition}
Let $U$ be a topological space and $(V, d_V)$ a metric space. 
We say that a sequence of maps $f_m\colon U \to V$ converges uniformly to $f\colon U\to V$ if for any $\varepsilon>0$ there exists $N\in\N$ such that $d_V( f_m(p) , f(p) )<\varepsilon$ for any $p\in U$. 

We say that $f_m$ converges locally uniformly if for any point $p\in U$ there is an open neighbourhood $B(p)$ of $p$ such that the sequence of maps $f_m|_{B(p)}$ converges uniformly on $B(p)$.
\end{definition}

Let $(f_m)_m\in\Bir(\p^n_{\k})$ be a sequence of birational transformations and suppose that we find a Euclidean open set $\Omega\subset\p^n_{\k}(\k)$ such that each $f_m$ is regular on $\Omega$. We can shrink $\Omega$ to an open subset of ${\k}^n$ and consider the regular map $f_m|_{\Omega}\colon\Omega \to \p^n_{\k}$, whose image lies in $\p^n_\k(\k)$.

Let $(f_m)_m\subset\Bir(\p^n_{\k})$ be a sequence of birational transformations and $\Omega\subset\p^n_{\k}(\k)$. 
We say that $f_m$ converges (locally) uniformly to $f\in\Bir(\p^n_{\k})$  on $\Omega$ if 
\begin{enumerate}
\item $f_m|_{\Omega}$ is regular for all $m\in\N$, 
\item $f|_{\Omega}$ is regular 
\item $f_m|_{\Omega}$ converges (locally) uniformly to $f|_{\Omega}$ on $\Omega$.
\end{enumerate}

\bigskip

Note however that there exist sequences $(f_m)_m\subset\Bir(\p^n_{\k})$ that do not converge, but such that there is a non-empty open subset $\Omega\subset\p^n_\k(\k)$ and $f\in \Bir(\p^n_{\k})$
such that each $f_m|_\Omega$ is regular and $f_m|_\Omega$ converges locally uniformly to $f|_\Omega$. This is illustrated by the following example:

\begin{example}
 Consider the sequence 
 \[
\C^n\rightarrow \C^n,\quad (X_1,X_2,\dots,X_n)\mapsto (X_1+\frac{1}{m!}X_2^m, X_2,\dots,X_n)
 \] 
 of automorphisms of $\C^n$. On every bounded subset $\Omega\subset \C^n$ this sequence converges uniformly towards the identity map.
 
 However, the induced sequence $(f_m)_m\subset \Bir(\p^n_\C)$, which is given by 
 $$f_m([x_0:x_1:x_2:\dots:x_n])=[x_0^m: x_0^{m-1}x_1+\frac{1}{m!}x_2^m:x_0^{m-1}x_2:\dots:x_0^{m-1}x_n],$$
 does not converge with respect
 to the Euclidean topology on $\Bir(\p^n_\C)$ since its degree is not bounded (cf. Propostion~\ref{lem:properties_topology}\ref{pro_top:4}).
 
This example can be generalised to any local field $\k$ by taking the coefficient of $X_2^m$ to be any sequence $(a_m)_{m\in\N}$ in $\k$ converging quickly enough to zero.
\end{example}

 Moreover, it might happen that a sequence $f_m\in \Bir(\p^n_{\k})$ converges to an automorphism $f\in \Aut(\p^n_{\k})$ and each $f_m$ is regular on some open subset $\Omega\subseteq \p^n_{\k}(\k)$,
 but the sequence of maps $f_m|_\Omega:\Omega\rightarrow \p^n_{\k}(\k)$ does not even converge pointwise to $f|_\Omega$ on all of~$\Omega$ as can be seen by the following example:
 
 \begin{example}\label{ex:non_global_point_convergence}
 Let $f_m\in \Bir(\p^2_\C)$ be defined by 
 $$f_m([x_0:x_1:x_2])=[x_0^2: x_0x_1+\frac{1}{m}x_2^2: x_0 x_2].$$
 This sequence converges to $f=\mathrm{id}\in\Bir(\p^2_\C)$ and each $f_m$ is regular on $\Omega=\p^2_\C\setminus\{[0:1:0]\}$.
 However, we have $f_m([0:1:p_2])=[0:1:0]$ for each $p_2\neq 0$ and thus $f_m([0:1:p_2])$ does not converge to $[0:1:p_2]=f([0:1:p_2])$.

Again, we can generalise this example for dimension $n\geq3$ and over other local fields by taking the coefficient of $x_2$ to be any sequence $(a_m)_{m\in\N}$ in $\k$ converging to zero.
\end{example}

Nevertheless, the next proposition shows that any convergent sequence of birational transformations of $\p^n$ has a subsequence that converges locally uniformly on some open subset 
of~$\p^n$.

\begin{proposition}\label{prop:uniform_convergence1}
Let $\k$ be a local field. Let $(f_m)_{m\in\N}\subset\Bir(\p^n_{\k})$ be a sequence converging to an element $f\in \Bir(\p^n_{\k})$.
Then, after passing to a subsequence, there is a non-empty open subset $\Omega\subseteq \p^n_{\k}(\k)$ such that
\begin{enumerate}
 \item\label{uniform11} $f|_{\Omega}:\Omega\rightarrow \p^n_{\k}$ and $f_{m}|_\Omega:\Omega\rightarrow \p^n_{\k}$ are regular for $m\gg 0$, and 
 \item\label{uniform12} the sequence $f_{m}|_{\Omega}:\Omega\rightarrow \p^n_{\k}(\k)$ converges to $f|_\Omega$ with respect to the compact-open topology on the space of continuous maps
 $\Omega\rightarrow \p^n_{\k}(\k)$, or equivalently $f_{m}|_\Omega$ converges locally uniformly to $f|_\Omega$ (with respect to any metric on $\p^n_{\k}(\k)$ inducing the Euclidean topology).
\end{enumerate}
\end{proposition}

Before proving Proposition~\ref{prop:uniform_convergence1}, let us first consider the following example which illustrates that it is not generally possible to find a Zariski open subset $\Omega\subseteq \p^n_{\k}(\k)$ satisfying the requirements of Proposition~\ref{prop:uniform_convergence1}.

\begin{example}
 Let $f_m\in \Bir(\p^2_\C)$ be defined by 
 $$f_m([x_0:x_1:x_2])=[x_0^2: x_0x_1+\frac{1}{m}x_2^2: x_0 x_2]$$
 as in Example~\ref{ex:non_global_point_convergence}.
 Moreover, let $\varphi_m\in\Bir(\p^2_\C)$ be a sequence of automorphisms defined by 
 $$\varphi([x_0:x_1:x_2]) = [x_0-\frac{1}{m}x_1: x_1:x_2], $$
 and consider the sequence of compositions $f_m\circ \varphi_m$.
 This sequence $f_m\circ\varphi_m$ also converges to $\mathrm{id}\in\Bir(\p^2_\C)$ since the individual sequences are both converging to $\mathrm{id}$.
 Each $f_m$ contracts the line $\{[x_0:x_1:x_2]\,|\,x_0 = 0\}$ to the point $[0:1:0]$ and $\varphi_m$ is mapping the line $\{[x_0:x_1:x_2]\,|\,x_0 + \frac{1}{m} x_1 = 0\}$ onto the line $\{[x_0:x_1:x_2]\,|\,x_0  = 0\}$.
 Hence, the composition $f_m\circ \varphi_m$ contracts the line $\{[x_0:x_1:x_2]\,|\,x_0 + \frac{1}{m} x_1 = 0\}$ to the point $[0:1:0]$.
 Since the lines $\{[x_0:x_1:x_2]\,|\,x_0 + \frac{1}{m} x_1 = 0\}$ are all distinct for different $m\in \mathbb N$ and local uniform convergence implies in particular pointwise convergence, there does not exist a Zariski open subset $\Omega\subseteq \p^2_\C$ such that the restriction $(f_m\circ \varphi_m)|_\Omega$ converges locally uniformly (or even only pointwise) to $\mathrm{id}|_\Omega$. 

 Again, we can generalise this example over other local fields by taking the coefficient of $x_2$ to be any sequence $(a_m)_{m\in\N}$ in $\k$ converging to zero. The example extends in an obvious way for $n\geq3$. 
\end{example}

\begin{proof}[Proof of Proposition~\ref{prop:uniform_convergence1}]
First, we argue that it is enough to consider the case where $f=\id$. If $f\neq \id $, consider 
 the sequence $f'_m= f_mf^{-1}$ which converges to $f'= f f^{-1}=\id$. If $\Omega $ is now a subset of $\p^n_\k(\k)$ satisfying \ref{uniform11} and \ref{uniform12} for $f'_m$,
 then $\Omega'=f^{-1}(\Omega \cap \mathrm{Def}(f^{-1}))$ gives the desired subset for the sequence $f_m$, where $\mathrm{Def}(f^{-1})$ denotes the largest subset of
 $\p^n_\k(\k)$ on which  $f^{-1}$ is regular.
 Thus, let us assume $f=\id $ in the following.
 
 By Proposition~\ref{lem:properties_topology}\ref{pro_top:4} there is a positive integer $d$ such that $f_m\in \Bir(\p^n_{\k})_{\le d}$ for all $m\gg 0$.
 Consider now the map $\pi_d: H_d(\k)\rightarrow \Bir(\p^n_{\k})_{\le d}$ from \S\ref{ss:zariskitopology}. 
 Since the map $\pi_d: H_d(\k)\rightarrow \Bir(\p^n_\k)_{\le d}$ is surjective and proper by Proposition~\ref{lem:properties_topology}\ref{pro_top:1}, there is a sequence $(p_m)_m\subset H_d(\k)$ such that $\pi_d(p_m)=f_m$ and 
 which converges to an element $p\in H_d(\k)$ with $\pi_d(p)=\id$ after possibly passing to a subsequence.
 
 We have $H_d(\k)\subset W_d(\k)$, where $W_d(\k)$ (again cf. Definition~\ref{DefWHG}) can be identified with the projectivisation $\p(V_d^{n+1}(\k))=\p(V_d(\k)\times \ldots\times V_d(\k))$ and $V_d=\k[x_0,\ldots, x_n]_d$ denotes the space of homogeneous polynomials of degree~$d$.
 
 Each $p_m$ is naturally an element of $W_d(\k)$ and we may find representatives $F_m=(F_{m,0},\ldots, F_{m,n})\in V_d(\k)\times\ldots\times V_d(\k)$, 
 each $F_{m,j}$ a homogeneous polynomial of degree $d$, such that we have
 $[F_{m,0}:\ldots :F_{m,n}]=p_m$ in $W_d(\k)$ and thus $[F_{m,0}:\ldots :F_{m,n}]=f_m$ as elements of $\Bir(\p^n_{\k})$.
 Furthermore, we can choose the $F_m$ in such a way that
e.g.\ the sum of the squares of all coefficients of $F_{m,j}$ for $j= 0,\ldots, n$ is equal to $1$. 
Then, after possibly again passing to a subsequence, 
 the sequence $F_{m}$  converges to an element $G=(G_0,\ldots,  G_n)\in V_d(\k)\times \ldots\times V_d(\k)$ (with respect to the Euclidean topology on the finite-dimensional vector space
 $V_d(\k)\times \ldots \times V_d(\k)$).
The sum of the squares of all coefficients of $G_j$ for $j = 0,\ldots,  n$ is also equal to~$1$, 
 hence the limit $G$ cannot be $0$ and $G =(G_0,\ldots, G_n)$ necessarily satisfies $[G_0: \ldots :G_n]=\id $ in $\Bir(\p^n_{\k})$.
 
After possibly again passing to a subsequence, the sequence of maps
 $$F_m:\k^{n+1}\rightarrow \k^{n+1},\, x=(x_0,\ldots, x_n)\mapsto F_m(x)=(F_{m,0}(x),\ldots, F_{m,n}(x))$$
 converges locally uniformly to the map 
 $$G:\k^{n+1}\rightarrow \k^{n+1},\, x=(x_0,\ldots, x_n)\mapsto G(x)=(G_0(x),\ldots, G_n(x)).$$
 Moreover, the equality $[G_0:\ldots:G_n]=\id$ in $\Bir(\p^n_{\k})$ implies that there is some homogeneous polynomial $H$ of degree $d-1$ by \cite[Lemma 2.13]{BF13}
 with $G=(G_0,\ldots, G_n)= (H x_0,\ldots,H  x_n)$ and thus $G(x)=H(x)x$.
 
Consider now the open subset $U=\k^{n+1}\setminus \{x\,|\, H(x)=0\}\subseteq \k^{n+1}\setminus \{0\}$. 
We have $H(G(x))=H(H(x)x)=H(x)^{d-1}H(x)=H(x)^{d}$ since $H$ is homogeneous of
 degree~$d-1$. This implies that $H(G(x))=0$ holds if and only if $H(x)=0$ holds and hence we have $G(U)=U$.
 Since $F_{m}$ converges locally uniformly to $G$, there is a non-empty open subset $U'\subseteq U$ such that $F_{m}(U')\subseteq U$ for all $m\gg 0$ and in particular
 $F_{m}(x)\neq 0$ for every $x\in U'$.
 
 Next, define $\Omega\subseteq \p^n_{\k}(\k)$ to be the image of $U'\subseteq \k^{n+1}\setminus \{0\}$ under the canonical projection. The subset $\Omega$ is open and non-empty,
 and by construction the birational transformations $f_m:\p^n_{\k}\rat \p^n_{\k}$, $[x_0,\ldots, x_n]\rat [F_{m,0}:\ldots :F_{m,n}]$ are regular when restricted 
 to $\Omega$ since $F_m(x)\neq 0$ for $x\in U'$.
 Moreover, as $F_m$ and hence $F_m|_{U'}$ converge locally uniformly to $G=(G_0,\ldots, G_n)=(H x_0,\ldots, H x_n)$, 
 the sequence $f_m|_\Omega:\Omega\rightarrow \p^n_{\k}(\k)$ converges locally uniformly to $f|_\Omega=\id|_\Omega:\Omega\rightarrow \p^n_{\k}(\k)$. 
\end{proof}

\begin{remark}\label{rmk:convergence1}
We extract the following from the proof of Proposition~\ref{prop:uniform_convergence1} and use it in the case where $\k$ is non-Archimedean and $f_m$ converges to $\id$. There exists $d\geq1$, a sequence $(F_m)_m$ in $(V_d(\k))^{n+1}$ and $G\in V_d(\k)^{n+1}$ such that $(f_m)_m$ and $\id$ are respectively the rational maps on $\p^n_\k$ induced by $F_m$ and $G$, and such that $F_m$ converges to $G$ in the Euclidean vector space $V_d(\k)^{n+1}$. There are moreover open sets $U'\subset U$ of $\k^{n+1}\setminus\{0\}$ such that $\p(U')=\Omega$, $F_m(U')\subset U$, $G(U)=U$ and $F_m(x)\neq0$ for all $x\in U'$ for all $m\gg0$.
\end{remark}

\subsection{Sequences of bounded order over $\C$}

\begin{proposition}\label{prop:sequences_of_bounded_order}
Let $\Omega\subset \C^n$ be an open connected subset and $f_m:\Omega\rightarrow \C^n$ a sequence of holomorphic maps locally uniformly converging to $\id|_\Omega:\Omega\rightarrow \C^n$.
Assume moreover that the sequence of holomorphic maps $f_m$ is of order $D>0$ in the sense that there is an open non-empty relatively compact subset $\Omega'\subset \Omega$
such that $(f_m)^i(\bar\Omega')\subseteq \Omega$ for all $i$-th iteraties of $f_m$, $i=1,\ldots, D$, and such that $(f_m)^{D}=\id$ on $\Omega'$ for all $m$.

Then $f_m=\id|_\Omega$ for $m\gg 0$.
\end{proposition}

The proof of this proposition will follow from the next two technical lemmas, the first of which makes use of Cartan's uniqueness theorem for holomorphic functions and the second of which uses Brouwer's fixed point theorem. The latter fails over $p$-adic numbers (for instance the map $f(x)=x+1$ sends the compact closed unit ball $\Z_p$ to itself but does not have fixed points on it).

\begin{lemma}\label{lemma:invariant_subset}
Let $\Omega$ and $\Omega'$ be open subsets of $\C^n$ and $f_m:\Omega\rightarrow \C^n$ a sequence of holomorphic maps,
all satisfying the same assumptions as in Proposition~\ref{prop:sequences_of_bounded_order}.
 Then for each $m\gg 0$ there is a non-empty compact convex subset $B_m\subset\Omega'$ with $f_m(B_m)=B_m$.
\end{lemma}

\begin{proof}
 Without loss of generality we may assume that $\Omega'$ is an open ball around $0\in \C^n$ of radius $r'$. Since $f_m|_\Omega$ converges locally uniformly to $\id|_\Omega$
 and all the iterates $f_m^i$ are defined on $\Omega'$, $f_m^i|_{\Omega'}$ converges locally uniformly to $\id|_{\Omega'}$ for each $i=1,\ldots, D-1$.
 
 Let $h_{m,i}:\Omega'\rightarrow \R$, $h_{m,i}(z)=\|f_m^i(z)\|^2$, where $\|\cdot \|$ denotes the Euclidean norm on $\Omega\subset \C^n$.
 Now fix some $r>0$ with $r<r'$ and set
 $$B_m=\bigcap_{i=1}^D h_{m,i}^{-1}([0,r^2])=\{z\in \Omega'\,|\, \|f_m^i(z)\|\leq r \text{ for } i=1,\ldots, D\}.$$
 First note that $B_m$ is contained in $\Omega'$ since $f_m^D=\id $ and $r<r'$.
 We have $0\in \Omega'$ and for each $i$ the sequence $f_m^i(0)$ converges to $\id(0)=0$ and thus $ B_m$ is non-empty for $m\gg 0$. 
 The compactness of $B_m$ follows from the fact that it is a closed subset of the closed ball of radius $r$ around $0\in \C^n$.
 Furthermore, we have $h_{m,i}(f_m(z))=h_{m,i+1}(z)$ for all $z\in B_m$ by construction, which implies 
 $$B_m\supseteq f_m(B_m)\supseteq f_m^2(B_m)\supseteq \ldots \supseteq f_m^D(B_m)=B_m$$
 and hence $f_m(B_m)=B_m$.
 It remains to show that $B_m$ is convex for $m\gg 0$.
 
 Since $f_m^i$ converges locally uniformly to $\id $ on $\Omega'$ for each $i$, all its partial derivates also converge locally uniformly to the partial derivatives of 
 $f_m^i$ and hence the sequences $(h_{m,i})_m$ converge locally uniformly to the function $h:\Omega'\rightarrow \R$, $h(z)=\|z\|^2$, and the partial derivatives of 
 $h_{m,i}$ converge to the corresponding partial derivatives of~$h$ for each $i=1,\ldots, D-1$.
 In particular, the Hessian of $h_{m,i}$ converges locally uniformly to the Hessian of $h$. As the latter is positive definite on $\bar\Omega$, 
 it follows that the Hessian of $h_{m,i}$ is positive definite for $m\gg0$. 
Consequently, $h_{m,i}$ is a convex function for $m\gg0$. In particular, $h_{m,i}^{-1}([0,r^2])$ is convex for $m\gg0$ and thus each
\[
h_{m,i}^{-1}([0,r^2])=\{z\in\Omega' \mid ||f_m^i(z)||\leq r\}
\]
 is convex. 
 Therefore, also the intersection 
\[
B_m=\bigcap_{i=1}^{D} h_{m,i}^{-1}([0,r^2])
\]
is convex for $m\gg 0$. 
\end{proof}

\begin{lemma}\label{lemma:fixed_point}
 Let $\Omega$ and $\Omega'$ be open subsets in $\C^n$, and $f_m:\Omega\rightarrow \C^n$ a sequence of holomorphic fuctions, all satisfying the same 
 assumptions as in Proposition~\ref{prop:sequences_of_bounded_order}. Let $B_m\subset\Omega'$ be a compact convex subset of $\Omega'$ with 
 $f_m(B_m)=B_m$ for $m\gg 0$, which exists by Lemma~\ref{lemma:invariant_subset}.
 
 Then for any $m\gg 0$, there is a fixed point $P_m=f_m(P_m)\in B_m$ of $f_m$, 
 and we have $Df_m(P_m)=\id$ for $m\gg 0$.
\end{lemma}

\begin{proof}
 The map $f_m|_{B_m}:B_m\rightarrow B_m$ is a continuous self-map of the compact convex set $B_m$ and thus has a fixed point $P_m\in B_m$ by Brouwer's fixed point theorem.
 
 Suppose that there are arbitrarily large $m$ such that $Df_m(P_m)\neq \id$, then after passing to a subsequence we may assume that 
 $Df_m(P_m)\neq \id $ for all $m$. Hence it is enough to prove the statement for a subsequence of $f_m$.
 
 Since $B_m\subset \bar\Omega'$ and $\bar\Omega'$ is compact, the sequence $P_m$ converges (after passing to a subsequence) to a point $P\in \bar\Omega'$.
 Consider now the sequence $Df_m(P_m)$ of differentials which can be canonically identified with a sequence of $n\times n$-matrices because $\Omega'\subset\C^n$.
 The identity $f_m^{D}=\id $ on $\Omega'$ implies $(Df_m(P_m))^D=(Df_m^{D})(P_m)=\id$.
 Consequently, all eigenvalues of $Df_m(P_m)$ have to be $D$-th roots of unity and $Df_m(P_m)$ is diagonisable.
 Since $f_m|_{\Omega'}$ converges locally uniformly to $\id $,
 also the differential $Df_m$ converges locally uniformly to $\id$ and hence
 the sequence of matrices $Df_m(P_m)$ converges to~$\id$. Therefore, all eigenvalues of $Df_m(P_m)$ are equal to $1$ for $m\gg 0$ and 
 hence $Df_m(P_m)=\id$ for $m\gg 0$.
\end{proof}

Using these two lemmas, we can now prove Proposition~\ref{prop:sequences_of_bounded_order}.

\begin{proof}[Proof of Proposition~\ref{prop:sequences_of_bounded_order}]
 By Lemma~\ref{lemma:invariant_subset} there are compact convex subsets $B_m\subset \bar\Omega'$ with $f_m(B_m)$ for $m\gg 0$
 and fixed points $P_m\in B_m$ of $f_m$ with $Df_m(P_m)=\id$ for $m\gg 0$.
 Now Cartan's uniqueness theorem (see e.g.\ \cite[p. 66]{N71}) 
 implies that $f_m=\id $ on the interior of $B_m$ and hence $f_m=\id$ on all of $\Omega$ by the identity principle.
\end{proof}

\subsection{Sequences of bounded order over non-Archimedean fields of characteristic zero}\label{ss:padic}

Let $\k$ be a non-Archimedean local field.
Let us denote, only in this subsection, by $|\cdot|$ the absolute value on $\k$. We denote by $\bR$ the closed unit ball in $\k$, and recall that $\bR$ is compact and also open and that it is the ring of integers of $\k$.

The algebra $\bR[x_1,\dots,x_n]$ carries the norm $||\sum a_Ix^I ||=\sup_I|a_I|$. Its elements map $\bR^n$ onto $\bR$ because $\bR$ is a ring.

The {\em Tate algebra} $\bR\langle x_1,\dots,x_n\rangle$ is the completion of $\bR[x_1,\dots,x_n]$ with respect to this norm. It is the set of all power-series over $\bR$ converging absolutely on $\bR^n$ and it is equal to the set of all elements $\sum a_Ix^I\in\bR[[x_1,\dots,x_n]]$ such that $|a_I|\to0$ if $I\to\infty$. In particular, the elements of $\bR\langle x_1,\dots,x_n\rangle$ map $\bR^n$ onto $\bR$.

The space $\k^N$ carries the norm $||(v_1,\dots,v_N)||=\max_i|v_i|$. 
Let $V_d(\k)\subset \bR[x_1,\dots,x_n]$ be the vector space of homogeneous polynomials of degree $d$ in $x_1,\dots,x_n$. If $V_d(\k)\simeq \k^N$ is the isomorphism given by a basis of monomials, then the norm on $V_d(\k)$ induced by this isomorphism is equal to the norm induced by $\bR[x_1,\dots,x_n]$.

\begin{lemma}\label{lem:padic1}
Let $\k$ be a non-Archimedean local field and let $f_m\in\Bir(\p^n_\k)$ be a sequence converging to $\id\in\Bir(\p^2_\k)$.  
Then there exists $\alpha\in\Aut(\p^n_\k)$ such that $f_m':=\alpha\circ f_m\circ\alpha^{-1}$ is regular on a Euclidean open subset of the chart $x_0\neq0$ for $m\gg0$. We write $f_m'=(f_{m,1}', \dots,f_{m,n}')$ where $f_{m,i}'$ are quotients of polynomials in $x_1,\dots,x_n$. Then we can moreover choose $\alpha$ such that $f_{m,i}'\in\bR\langle x_1,\dots,x_n\rangle$ and $f_{m,i}'$ converges to $x_i$ in $\bR\langle x_1,\dots,x_n\rangle$ for $i=1,\dots,n$.
\end{lemma}
\begin{proof}
By Proposition~\ref{prop:uniform_convergence1} and Remark~\ref{rmk:convergence1}, there exists $d\geq1$ and a sequence $F_m$ in $V_d(\k)^{n+1}$ and $G\in V_d(\k)^{n+1}$ such that their induced maps on $\p^n_\k$ are $f_m$ and $\id$, respectively, and such that $F_m$ converges to $G$ in the Euclidean topology. Moreover, there exist Euclidean open subsets $U'\subset U\subset\k^{n+1}\setminus\{0\}$ such that $F_m(U')\subset U$, $G(U)= U$ and $F_m(x)\neq0$ for all $x\in U'$ and $m\gg0$ and such that $f_m$ is regular on $\Omega=\p(U')$ for $m\gg0$.

There is an affine transformation $\alpha$ of $V_d(\k)^{n+1}$ such that $(1,\dots,1)\in\alpha(U')$. 
Then $F_m':=\alpha\circ F_m\circ\alpha^{-1}=[F_{m,0}'(x_0,\dots,x_n):\dots: F_{m,n}'(x_0,\dots,x_n)]$ converges to $G':=\alpha\circ G\circ\alpha^{-1}$ in $V_d(\k)^{n+1}$. 
Let $f_m'$ be the rational map of $\p^n$ induced by $F_m'$. Note that $G'$ induces the identity on $\p^n_\k$.
We have $F_m'(x)\neq0$ for any $x\in\alpha(U')$ and hence $f_m'$ is birational. It is moreover regular on $\p(\alpha(U'))=\alpha(\Omega)\ni [1:\dots:1]$ for $m\gg0$ and in particular, it is regular on a Euclidean open set in $\p^n\setminus\{x_0=0\}$ for $m\gg0$.

Let us show that for this choice of $\alpha$, the maps $f_m'=(f_{m,1}',\dots,f_{m,n}')$ are contained in $\bR\langle x_1,\dots,x_n\rangle$ for $m\gg0$.
There exists $\lambda\in\k^*$ such that $\lambda G'=(\lambda G_0',\dots,\lambda G_n')\in\bR[x_0,\dots,x_n]^{n+1}$. Then the sequence $(\lambda F_m')_m$ converges to $\lambda G'$ and $\lambda F_m'$ induces the birational map $f_m'$ for each $m$. Since $\bR$ is open in $\k$, we have $\lambda F_m'\in\bR[x_0,\dots,x_n]^{n+1}$ for $m\gg0$. Therefore, on the chart $x_0=1$ and for $m\gg0$, the quotients $f_{m,i}'=\frac{\lambda F_{m,i}'(1,x_1,\dots,x_n)}{\lambda F_{m,0}'(1,x_1,\dots,x_n)}$ are quotients of polynomials in $\bR[x_1,\dots,x_n]$.
In particular, they are power-series contained in $\bR[[x_1,\dots,x_n]]$. 
We write $f_{m,i}'=\sum_I a_{I,i}\ x^I$. 
By hypothesis on $\alpha$, the map $f_m'$ is regular at the point $(1,\dots,1)\in\alpha(U')$.  
Since $f_{m,i}'(1,\dots,1)=\sum a_{I,i}$, it follows that $|a_{I,i}|$ converges to zero when $I\to\infty$. This means that $f_{m,i}'$ is contained in the Tate algebra $\bR\langle x_1,\dots,x_n\rangle$ for $m\gg0$. 

It remains to show that $f_{m,i}'$ converges to $x_i$ in $\bR\langle x_1,\dots,x_n\rangle$ for $i=1,\dots,n$. 
By hypothesis, for each $i$ the sequence $F_{m,i}'$ converges to $G_i'$ in $V_d(\k)$ and hence the sequence $(\lambda F_{m,i}')_m$ converges to $\lambda G_i'$ in $\bR\langle x_0,\dots,x_n\rangle$.
Therefore, the sequence $\lambda F_{m,i}'(1,x_1,\dots,x_n)$ converges to $\lambda G_i'(1,x_1,\dots,x_n)$ in $\bR\langle x_1,\dots,x_n\rangle$. 
The quotient  
$f_{m,i}'(x_1,\dots,x_n)=\frac{\lambda F_{m,i}'(1,x_1,\dots,x_n)}{\lambda F_{m,0}'(1,x_1,\dots,x_n)}$ is contained in $\bR\langle x_1,\dots,x_n\rangle$ and therefore converges to $\frac{\lambda G_i'(1,x_1,\dots,x_n)}{\lambda G_0'(1,x_1,\dots,x_n)}=x_i$. 
\end{proof}

\subsection{Proof of Theorem~\ref{thm1}}

Any local field of characteristic zero is either isomorphic to $\R$, to $\C$ or to a finite extension of the $p$-adic numbers $\Q_p$, see for instance \cite[Remark 7.49, p.127]{M20}. The proof of Theorem~\ref{thm1} for finite extensions of $\Q_p$ is based on a corollary of the Bell-Poonen theorem \cite{P14} that is proven in \cite{CX18}, and which holds over characteristic zero. 

\begin{proof}[Proof of Theorem~\ref{thm1}]
Let $\k$ be a local field of characteristic zero. 
Since $\ord(f_m)\leq D$ for all $m$, we have $(f_m)^N=\id$ for all $m$ and $N=D!$.
By Proposition~\ref{prop:uniform_convergence1}\ref{uniform11}\&\ref{uniform12} and after possibly passing to a subsequence, there is an open subset $\Omega_0\subset \p^n_\k(\k)$ such that 
$f_m|_{\Omega_0}$ is regular for all $m\gg 0$ and $f_m|_{\Omega_0}:\Omega\rightarrow \p^n_\k(\k)$ converges locally uniformly to $\id|_{\Omega_0}$.

{\bf Suppose that $\k=\R$ or $\k=\C$}. By Lemma~\ref{lem:subset-topology} it suffices to prove the claim for $\k=\C$.
After shrinking we may assume that $\Omega_0$ is biholomorphic to an open ball around $0$ in $\C^n$ with radius $r_0>0$ and we fix such an identification $\Omega_0\subset \C^n$.

For any $r>0$ with $r<r_0$, let $\Omega_r\subset \Omega_0$ be the open ball around $0$ with radius~$r$.
Let $r_j>0$, $j=1,\ldots, N$, be positive numbers with 
$r_0>r_1>\ldots> r_N$. These give rise to a strictly increasing sequence of open sets
$\Omega_{r_N}\subset \Omega_{r_{N-1}}\subset \ldots \subset \Omega_{r_1}\subset\Omega_{r_0}=\Omega_0$, and the closure $\bar\Omega_{r_j}$ is contained in $\Omega_{r_{j-1}}$ for each 
$j=1,\ldots, N$. 
Since the sequence of maps $f_m|_{\Omega_0}:\Omega_0\rightarrow \p^n$ converges locally uniformly to $\id|_{\Omega_0}:\Omega_0\rightarrow \p^n$, 
we have $f_m(\bar\Omega_{r_j})\subset \Omega_{r_{j-1}}$ for each $j=1,\ldots, N$ and $m\gg 0$. 
This implies in particular that 
$f_m^i|_{\Omega_{r_N}}$ is a regular map $\Omega_{r_N}\rightarrow \Omega_0\subset \C^n$ for $m\gg 0$ and each $i=1,\ldots , N$. And by assumption we also have 
$f_m^N|_{\Omega_{r_N}}=\id|_{\Omega_{r_N}}$. 
An application of Proposition~\ref{prop:sequences_of_bounded_order} now yields $f_m|_{\Omega_0}=\id $ for $m\gg 0$ and hence $f_m=\id$ on $\p^n$ by the identity principle.

{\bf Suppose that $\k$ is a finite extension of $\Q_p$}.  
Lemma~\ref{lem:padic1} applied to our sequence $f_m$ implies that we can choose $\alpha\in\Aut(\p^n_\k)$ such that $f_m':=\alpha\circ f_m\circ\alpha^{-1}$ 
is regular on an open subset of the chart $x_0\neq0$ for $m\gg0$. 
We write $f_m'$ as map $f_m'=(f_{m,1}', \dots,f_{m,n}')$ given by quotients $f_{m,i}'$ of polynomials in $x_1,\dots,x_n$. 
Lemma~\ref{lem:padic1} states moreover that we can choose $\alpha$ such that $f_{m,i}'\in\bR\langle x_1,\dots,x_n\rangle$ for $m\gg0$ and such that $f_{m,i}'$ converges to $x_i$ for $i=1,\dots,n$.

In particular, the sequence $||f_{m,i}'-x_i||$ is a Cauchy sequence in $\bR\langle x_1,\dots,x_n\rangle$ for $i=1,\dots,n$. 
If we set $\varepsilon=\frac{1}{p^2}$, then then there exists $m_0\in\N$ such that $||f_{m,i}'-x_i||<1/p^2$ for $m\geq m_0$. 
The norm on $\bR\langle x_1,\dots,x_n\rangle^n$ is given by $||(h_1,\dots,h_n)||=\max_i||h_i||$. It follows that
\begin{equation}\tag{$\ast$}\label{ast}
||f_m'-\id||\leq 1/p^2,\quad m\gg0
\end{equation}
\cite[Corollary 2.5]{CX18} (with $c=2>1/(p-1)$) implies $f_m'=\id$ for $m\gg0$  (here we are using $\mathrm{char}(\k)=0$). It follows that $f_m=\id$ for $m\gg0$. 
\end{proof}
Condition (\ref{ast}) is also referred to as $f_{m}\equiv \id$ (mod $p^2$).

\begin{proof}[Proof of Corollary~\ref{thm2}]
A morphism $\rho\colon\A^1_{\k}\to\Bir(\p^n)$ is continuous in the Euclidean topology \cite[Lemma 2.11]{BZ16}.
Hence the sequence $\rho(t_m)$ converges to $\id$ in the Euclidean topology and by Theorem~\ref{thm1} we have $\rho(t_{m})=\id$ for $m\gg 0$.
The morphism $\rho$ is also continuous in the Zariski topology (by its definition), so $\rho^{-1}(\id)\subset\A^1_{\k}$ is a Zariski closed subset containing infinitely many points.
If follows that $\rho^{-1}(\id)=\A^1_{\k}$.
\end{proof}


\section{Sequences of converging subgroups}
Let $\k$ be  local field. 
In this paragraph, we study properties of sequences of subgroups in $\Bir(\p^n_{\k})$ which converge to the identity $\id$ 
in the following sense:

\begin{definition}
 Let $(G_m)_{m\in \N}$ be a sequence of subgroups of a topological group~$G$.
 We say that $(G_m)_{m\in \N}$ converges to the neutral element $\id\in G$ if
 for any open neighbourhood~$U$ of $\id$ there exists $m_0\in \N$ such that $G_m \subset U$ for all $m> m_0$.
\end{definition}

If $G$ is a topological group containing a sequence $(G_m)_{m\in\N}$ of subgroups that converges to $\id\in G$ and such that $G_m\neq\{\id\}$ for each $m\geq1$ , then $G$ has small subgroups.

\begin{remark}\label{rmk: sequences in converging subgroups}
If $(G_m)_{m\in \N}$ is a sequence of groups converging to the neutral element $\id$ and $(g_m)_{m\in\N}$ is a sequence of elements with $g_m\in G_m$,
then the sequence $(g_m)_{m\in \N}$ converges to the element $\id$ in the usual sense.
\end{remark}

First, we show that a converging sequence of subgroups of $\Bir(\p^n_{\k})$ has to be bounded eventually.

\begin{lemma}\label{lem:subgroup_sequence_bounded_degree}
Consider a sequence of subgroups $(G_m)_m\subset \Bir(\p^n_{\k})$ converging to $\id\in\Bir(\p^n_{\k})$.
 Then for $m\gg 0$ the groups $G_m$ have simultaneously bounded degree in the sense that there is $m_0>0$ and $d>0$ such that $G_m\subset \Bir(\p^n_{\k})_{\leq d}$ for all $m\geq m_0$.
\end{lemma}

\begin{proof}
 Assume the statement of the lemma was wrong. Then there is a sequence of elements $(g_m)_{m\in \N}$ with $g_{m}\in G_{m}$ for all $m$ 
 such that the sequence of degrees $\deg(g_{m})$ is unbounded. By Proposition~\ref{lem:properties_topology}\ref{pro_top:4} this implies 
 that the sequence $(g_m)_{m\in \N}$ does not converge in $\Bir(\p^n_{\k})$, whereas by 
 Remark~\ref{rmk: sequences in converging subgroups} 
 $(g_m)_{m\in \N}$ converges to $g\in\Bir(\p^n_{\k})$, yielding a contradiction.
\end{proof}

If $G\subset \Bir(\p^n_{\k})$ is a subgroup which is closed with respect to the Zarisiki topology and of bounded degree, 
i.e.\ $G\subset \Bir(\p^n_{\k})_{\leq d}$ for some $d$, then 
$G$ carries the structure of an algebraic group compatible with the structure of $\Bir(\p^n_{\k})$ as proven in \cite[Corollary~2.18]{BF13}.
Analogously, we prove that over $\k=\R$ and $\k=\C$, subgroups of $\Bir(\p^n_{\k})$ which are closed with respect to the Euclidean topology of $\Bir(\p^n_{\k})$ and of bounded degree carry the structure of a {\em real} Lie group. 

\begin{proposition}\label{prop:closed_bounded_subgroups}
Let $\k=\C$ or $\k=\R$. 
Let $G\subset \Bir(\p^n_{\k})$ be a subgroup which is closed with respect to the Euclidean topology of $\Bir(\p^n_{\k})$ and of bounded degree in the sense that there is $d>0$ such that
 $G\subset \Bir(\p^n_{\k})_{\leq d}$. Then $G$ carries a unique structure of a {\em real} Lie group which is compatible, 
 i.e.\ whose induced topology on $G$ coincides with the restriction of the Euclidean topology of $\Bir(\p^n_{\k})$ to $G$.
\end{proposition}

\begin{proof}
 If a Lie group structure on $G$ exists, then it is unique since topological groups admit at most one Lie group structure compatible with 
their topology. 
 
 Let $H$ denote the closure of $G$ in $\Bir(\p^n_{\k})$ with respect to the Zariski topology. The subset $H$ is a closed subgroup of $\Bir(\p^n_{\k})$, which follows from the fact that the map
 $\Bir(\p^n_{\k})\times\Bir(\p^n_{\k})\rightarrow \Bir(\p^n_{\k})$, $(\psi,\phi)\mapsto\psi\phi^{-1}$,
 is continuous with respect to the Zariski topology (cf.\ remark after \cite[Definition~2.2]{BF13}).
 Since we have $G\subset \Bir(\p^n_{\k})_{\leq d}$ for some $d$ and since
 $\Bir(\p^n_{\k})_{\leq d}$ is closed with respect to the Zariski topology by \cite[Corollary~2.8]{BF13}, we also have $H\subset \Bir(\p^n_{\k})_{\leq d}$.
 
Hence, $H$ carries a unique structure of a real resp.\ complex Lie group structure for $\k = \R, \C$. Since the induced structure of any topologically closed subgroup of a (real or complex) Lie group gives rise to the structure of a {\em real} Lie group on that subgroup, the group $G$ carries the structure of a {\em real} Lie group compatible with the restricted Euclidean topology of $\Bir(\p^n_{\k})$ to $G$.
\end{proof}

\begin{remark}
Let $\k=\C$ or $\k=\R$. 
 If $G\subset \Bir(\p^n_{\k})$ is a bounded, but not necessarily closed subgroup in the sense of the above proposition, we can take the topological closure $\overline{G}$ 
 of $G$ with respect to the Euclidean topology. The closure $\overline{G}$ is again a subgroup and satisfies the assumptions of the proposition.
 Therefore, the closure $\overline{G}$ of any subgroup $G\subset \Bir(\p^n_{\k})$ of bounded degree carries a unique compatible structure of a real Lie group.
\end{remark}

An easy consequence of the preceding proposition is the following statement, compare with \cite[Proposition 5.11]{BF13}.
It is an analogue to the fact that every closed subgroup of a Lie group is a Lie group itself.

\begin{corollary}\label{cor:compact_subgroups}
Let $\k=\C$ or $\k=\R$. Let $G\subset \Bir(\p^n_{\k})$ be a compact subgroup (with respect to the Euclidean topology). Then $G$ carries a compatible real Lie group structure. 
\end{corollary}

\begin{proof}
 This statement directly follows from the fact that compact subsets of $\Bir(\p^n_{\k})$ are bounded (cf.\ Proposition~\ref{lem:properties_topology}\ref{pro_top:4}) and the preceding Proposition~\ref{prop:closed_bounded_subgroups}.
\end{proof}

It is worth noting that if $\k$ is a perfect field, then a compact subgroup of $\Bir(\p^n_{\k})$ that carries a structure of algebraic subgroup of $\Bir(\p^n_{\k})$ is finite. Indeed, algebraic subgroups of $\Bir(\p^n_{\k})$ are affine; this is shown in \cite[Remark 2.21]{BF13} using Chevalley's structure theorem for algebraic groups, which holds over perfect fields. Hence algebraic subgroups of $\Bir(\p^n_{\k})$ are non-compact in their natural Euclidean topology unless they are finite.

We can now use these results to prove Theorem~\ref{theorem:converging subgroups}. 

\begin{proof}[Proof of Theorem~\ref{theorem:converging subgroups}]
Let $\k=\R$ or $\k=\C$. 
 Let $(G_m)_{m\in \N}$ be a sequence of subgroups of $\Bir(\p^n_{\k})$ that converges to $\id \in \Bir(\p^n_{\k})$. If there is a subsequence with all $G_m$ finite, we apply Theorem~\ref{thm1}. So, let us assume that the $G_m$ are non-finite for $m\gg0$. For each $m$ let $H_m=\overline{G}_m$ denote the topological closure of $G_m$ in
 $\Bir(\p^n_{\k})$. By Lemma~\ref{lem:subgroup_sequence_bounded_degree} there is $m_0>0$ and $d>0$ such that $G_m\subset \Bir(\p^n_\k)_{\leq d}$ for all $m\geq m_0$.
 Since $\Bir(\p^n_{\k})_{\leq d}$ is closed in $\Bir(\p^n_{\k})$ (cf.\ Proposition~\ref{lem:properties_topology}\ref{pro_top:3}), this implies that $H_m=\overline{G}_m\subset \Bir(\p^n_{\k})_{\leq d}$ 
 for $m\geq m_0$. Without loss of generality we can thus assume that $G_m, H_m\subseteq \Bir(\p^n_{\k})_{\leq d}$ for all $m\in\N$.
 The local compactness of $\Bir(\p^n_{\k})_{\leq d}$ (cf.\ Proposition~\ref{lem:properties_topology}\ref{pro_top:2}) yields that we may assume that each $H_m=\overline{G}_m$ is compact since~$G_m$ converges to $\id\in \Bir(\p^n_{\k})_{\leq d}\subset \Bir(\p^n_{\k})$. 
Since $\Bir(\p^n_{\k})_{\leq d}$ is locally compact, it follows that $(H_m)_{m\in \N}$ converges to $\id$ in $\Bir(\p^n_{\k})_{\leq d}$.  By the construction of the topology, the sequence $(H_m)_{m\in \N}$ converges to $\id$ in $\Bir(\p^n_{\k})$ if and only if $(H_m)_{m\in\N}$ converges to $\id$ in $\Bir(\p^n_{\k})_{\leq d}$, and hence $(H_m)_{m\in\N}$ converges to $\id$ in $\Bir(\p^n_\k)$. 
 Consequently, it is enough to prove the statement of the theorem in the case where each subgroup $G_m$ is compact, which we will assume in the following. 

 By the preceding Corollary~\ref{cor:compact_subgroups}, each compact $G_m$ thus carries the structure of a Lie group. Every compact real Lie group $G\neq\{\id\}$ contains 
 non-trivial elements of order $2$ (in a compact subtorus of $G$). Theorem~\ref{thm1} implies now that $G_m=\{\id\}$ for $m\gg 0$ because otherwise we could construct a sequence of elements $g_m\in \Bir(\p^n_{\k})$
of order $2$ and converging to $\id$, but without $g_m=\id$ for all sufficiently large $m$. 
\end{proof}

\begin{remark}\label{rmk:fail for p-adic}
The completion $\C_p$ of the algebraic closure of the $p$-adic numbers is isomorphic to $\C$ as a field, but not as local fields, because $\C_p$ is not locally compact. The isomorphism induces an isomorphism of groups $\Bir(\p^n_{\C_p})\simeq\Bir(\p^n_\C)$ that is not a homeomorphism in the Euclidean topology: for instance $\Bir(\p^n_{\C_p})$ contains $\mathrm{PSL}_{n+1}(\Z_p)$, which contains arbitrary small subgroups. An example is the sequence of groups 
 \[
 G_m=\left\{
 \left(\begin{matrix}
 1&q&0&\dots\\
 0&1&0&\dots\\
 \vdots&&\ddots&\\
 0&\dots&0&1
 \end{matrix}\right)\bigg| \pm
 q=\sum_{i=k}^\infty a_i  p^i,\ a_i\in\{0,\ldots, p-1\}, \ k\geq m\right\}
 \] 
 which converges to $\id$. 
 In particular, $\Bir(\p^n_{\Q_p})$ contains small subgroups and Theorem~\ref{theorem:converging subgroups} fails over the p-adic numbers. 
\end{remark}

\section{Obstructions to lifting continuous maps}

Proposition~\ref{prop:uniform_convergence1} raises the question, whether we can prove a similar result for continuous maps $[0,1]\to\Bir(\p^n_{\k})$, that is, for continuous families of birational maps. The following example shows that it is not possible in general. 
The example is given for $n=2$ but can be adapted easily for any $n\geq3$.

\begin{example}\label{ex: oscillating map}
 Let $f_t\in\Bir(\p^2_\C)$, $t\in [0,1]$, be the birational tranformation given by 
 $$f_t\colon[x_0:x_1:x_2]\mapsrat[x_0^2:x_0x_1:x_0x_2+t(1-t)x_1^2].$$
 We have $f_0=f_1=\id$ and for $t\neq 0,1$ the transformations $f_t$ are not defined at $p=[0:0:1]$ and are regular everywhere else.
 
 Let $\sigma:[-1,1]\rightarrow \mathrm{PSU}(3)$ be a continuous surjective map. Such a space filling curve exists by the {\sc Hahn–Mazurkiewicz} theorem (see for instance \cite[Theorem 6.8]{Sagan}), because $\mathrm{PSU}(3)$ is compact, connected, locally connected and second countable. Using the inclusions $\mathrm{PSU}(3)\subset \mathrm{PGL}(3)\cong \Aut(\p^2_\C)\subset \Bir (\p^2_\C)$, we can view each $\sigma(t)$ as an element of 
 $\Bir(\p^2_\C)$. 
 Consider now the map $\hat\sigma:(0,1]\rightarrow \mathrm{PSU}(3)$, $\hat\sigma(t)=\sigma(\sin(1))^{-1}\cdot\sigma(\sin(\frac{1}{t}))$. 
 By construction $\hat\sigma|_{(0,\varepsilon)}:(0,\varepsilon)\rightarrow\mathrm{PSU}(3) $ is surjective
 for any $\varepsilon >0$.
 
 Consider now the continuous map $\rho:(0,1]\rightarrow \Bir(\p^2_\C)$, 
 $$\rho(t)=\hat\sigma(t)\circ f_t \circ (\hat\sigma(t))^{-1}
 =\sigma(\sin(1))^{-1}\circ\sigma\big(\sin\left(\frac{1}{t}\right)\big)\circ f_t\circ \big(\sigma\big(\sin\left(\frac{1}{t}\right)\big)\big)^{-1}\circ\sigma(\sin(1)).$$
 \end{example}
 
 \begin{lemma}\label{lem:counter_example1}
 The map $\rho$ from Example~\ref{ex: oscillating map} extends to a continuous map $\hat{\rho}\colon [0,1]\to\Bir(\p^n_\C)$ with $\hat{\rho}(0)=\id$.   
Moreover, the map 
\[H\colon[0,1]\times[0,1]\to\Bir(\p^n_\C),\quad 
H(s,t)=
\begin{cases} 
\hat{\rho}(t),& t\geq s\\ 
\hat{\sigma}(s)\circ f_t\circ \hat{\sigma}(s)^{-1},& t< s
\end{cases}
\]
is continuous. In particular, the map $\hat{\rho}$ is homotopic to the map $t\mapsto f_t$.
 \end{lemma}
 
 \begin{proof}
 The action $\varphi:\mathrm{PSU}(3)\times \Bir(\p^2_\C)\rightarrow \Bir(\p^2_\C)$, 
 $\varphi(g,f)=g\circ f\circ g^{-1}$, on $\Bir(\p^2_\C)$ is continuous and $f_0=\id$ is a fixed point of this action.
 Let $\Omega \subset \Bir(\p^2_\C)$ be any open neighbourhood of $f_0=\id$ in $\Bir(\p^2_\C)$.
 Then for any $g\in \mathrm{PSU}(3)$ there exists an open neighbourhood $U_g\subseteq \mathrm{PSU}(3)$ of $g\in \mathrm{PSU}(3)$ and an open 
 neighbourhood $V_g\subseteq \Bir(\p^2_\C)$ of $f_0$ such that $\varphi(U_g\times V_g)\subseteq \Omega$.
 Since $\mathrm{PSU}(3)$ is compact, there are finitely many $g$, say $g_1,\ldots, g_r$, such that 
 $\bigcup_{j=1} U_{g_j}=\mathrm{PSU}(3)$. Define now $V=\bigcap_{j=1}^r V_{g_j}$ and note that $V\subset V_{g_j}$ for all $j=1,\ldots, r$.
 Now we have by construction
 $$\varphi(\mathrm{PSU}(3)\times V)=\bigcup_{j=1}^r \varphi(U_{g_j}\times V)\subseteq 
 \bigcup_{j=1}^r \varphi(U_{g_j}\times V_{g_j})\subseteq \Omega.$$
 Since $[0,1]\rightarrow \Bir(\p^2_\C)$, $f\mapsto f_t$, is continuous, there is $\delta>0$ such that 
 $f_t\in V$ for all $t\in [0,\delta)$,  
 thus also $\hat{\rho}(0)=\id\in \varphi(\mathrm{PSU}(3)\times V)\subset\Omega$ and 
 $$\rho(t)=\hat\sigma(t)\circ f_t \circ (\hat\sigma(t))^{-1}=\varphi(\hat\sigma(t), f_t)\in \varphi(\mathrm{PSU}(3)\times V)\subseteq \Omega,\quad t\in(0,\delta),$$
 which shows that $\hat{\rho}^{-1}(\Omega)$ is a neighbourhood of $t=0$. This proves the first claim.
 
Let us show that $H$ it is continuous. 
Since $\hat{\rho}$ is continuous, it follows that $H$ is continuous on $\{(s,t)\in[0,1]\times[0,1]\mid t\neq0,t> s\}$. 
Composition and inversion in $\Bir(\p^2_\C)$ are continuous in the Euclidean topology (see Proposition~\ref{lem:properties_topology}\ref{pro_top:5}) and $[0,1]\mapsto\Bir(\p^2_\C)$, $f\mapsto f_t$ is continuous, and therefore $H(s,t)=\hat{\sigma}(s)\circ f_t\circ\hat{\sigma}(s)^{-1}$ is continuous on $\{(s,t)\in[0,1]\times[0,1]\mid s\neq0,t<s\}$. Moreover, for any $t\neq0$ we have $\lim_{s\rightarrow t}H(s,t)=H(t,t)$.  
It remains to check that $H$ is continuous in $(0,0)$.
We have
 \[
H(s,t)=
\begin{cases}
\hat{\rho}(t)=\varphi(\hat\sigma(t), f_t), & t\geq s,\ t\in[0,\delta)\\
\hat{\sigma}(s)\circ f_t\circ\hat{\sigma}(s)^{-1}=\varphi(\hat\sigma(s), f_t),& t<s, \ t\in[0,\delta)
\end{cases}
 \]
and hence $H(s,t)\in \varphi(\mathrm{PSU}(3)\times V)\subseteq \Omega$ for all $[0,1]\times[0,\delta)$. In particular, $H$ is continuous in $(s,t)=(0,0)$.

Finally, to conclude that $\hat{\rho}$ is homotopic to $t\mapsto f_t$, it suffices to see that $H(0,t)=\hat{\rho}(t)$ and $H(1,t)=\hat{\sigma}(1)\circ f_t\circ\hat{\sigma}(1)^{-1}=f_t$ for any $t\geq0$.
 \end{proof}
 
The second item of the following statement is Proposition~\ref{prop:example1-main}.

\begin{proposition}\label{prop:counter_example1}
Let $\hat{\rho}:[0,1]\to \Bir(\p^2_\C)$ be defined as in Lemma~\ref{lem:counter_example1}. Then the following hold.
\begin{enumerate}
\item
For any $q\in \p^2_\C$ and any $\varepsilon>0$ there exists a point $t_0\in (0,\varepsilon)$ such that the transformation $\hat{\rho}(t_0)$ is not regular at $q$.
\item
For any $n\geq2$ there does not exist an open non-empty subset $U\subset[0,1]\times\p^n_\C$ such that $U\to [0,1]\times\p^n_\C$, $(t,p)\mapsto(t,\hat{\rho}(t)(p))$ is a
well defined map. 
\end{enumerate}
\end{proposition}
\begin{proof}
(1)  Since $\mathrm{PSU}(3)$ acts transitively on $\p^2_\C$ and $\hat\sigma|_{(0,\varepsilon)}$ is surjective for any $\varepsilon>0$, there is some $t_0\in (0,\varepsilon)$ such that 
 $q=\hat\sigma(t_0)(p)$ for $p=[0:0:1]$. Because $f_{t_0}$ is not regular at $p$, it now follows that $\rho(t_0)=\hat\sigma(t_0)\circ f_{t_0}\circ (\hat\sigma(t_0))^{-1}$ 
 is not regular at $q=\hat\sigma(t_0)(p)$.
 
 (2) If $n\geq3$, we simply extend Example~\ref{ex: oscillating map} and Lemma~\ref{lem:counter_example1} to $\p^n_\C$ by using $\mathrm{PSU}(n)$. Then the claim follows directly as in (1). 
\end{proof}

\begin{remark}\label{rmk:Cantor}
The same example can be constructed over $\R$ by using $\mathrm{PSO}(3)$ instead of $\mathrm{PSU}(3)$, because $\mathrm{SO}(3)$ acts transitively on the $3$-sphere.

A similar example can be constructed for non-Archimedean local fields, whose valuation ring $\mathcal{O}$ is a Cantor set. In that case, $\mathrm{PGL}_3(\mathcal{O})$ is a compact Lie group that is a Cantor set. 
Let $C$ be the Cantor set in $[0,1]$ and let 
\[
C'=\{c\in (0,1]\mid \sin\left(\frac{1}{c}\right)\in C\cup -C\}\subset(0,1].
\] 
In Example~\ref{ex: oscillating map} we replace $\rho$ by the map $\rho\colon C'\to\mathrm{PGL}_3(\mathcal{O})$,
\[
\rho(c)=\sigma(\sin(1))^{-1}\circ\sigma\big(\sin\left(\frac{1}{c}\right)\big)\circ f_c\circ \big(\sigma\big(\sin\left(\frac{1}{c}\right)\big)\big)^{-1}\circ\sigma(\sin(1)).
\]
where $\sigma\colon C\cup -C\to C\stackrel{\sim}\to\mathrm{PGL}_3(\mathcal{O})$ is the composition of a continuous surjective map and a homeomorphism.  
\end{remark}

Proposition~\ref{prop:counter_example1} says that in general, continuous maps from compact sets to $\Bir(\p^n_\C)$ do not satisfy a condition similar to the definition of morphisms into $\Bir(\p^n_\C)$. 
However, since the map $\hat{\rho}$ from Lemma~\ref{prop:counter_example1} is homotopic to the map $t\mapsto f_t$ and all $f_t$ are regular on $\p^2_\C\setminus\{[0:0:1]\}$, we naturally arrive at the following questions. 

\begin{question}
Let $\rho\colon[0,1]\to\Bir(\p^n_\k)$ be a continuous map. Can we find a continuous map $\rho'\colon[0,1]\to\Bir(\p^n_\k)$ that is homotopic to $\rho$ and non-empty open subsets $U,V\subset[0,1]\times\p^n_\k$ such that $U\to V$, $(t,p)\mapsto(t,\rho'(t)(p))$ is a homeomorphism? 
\end{question}

\begin{question}
Let $\rho\colon[0,1]\to\Bir(\p^n_\k)$ be a continuous map. Can we find a continuous map $\rho'\colon[0,1]\to\Bir(\p^n_\k)$ that is homotopic to $\rho$, and a point $p\in\p^n_\k$ such that $\rho'(t)$ is regular at $p$ for any $t\in[0,1]$? 
\end{question}

It is unclear whether all continuous maps $[0,1]\to\Bir(\p^n_\k)$ are homotopic to an analytic path.

A much used property of any morphism $\rho$ from a variety $A$ to $\Bir(\p^n_\k)$ is that there is an open affine cover $A=\cup A_i$ and morphisms $\rho_i\colon A_i\to H_{d_i}$ such that $\rho|_{A_i}=\pi_{d_i}\circ\rho_i$.  This fails in general for continuous maps, as the next example shows.

\begin{example}\label{ex: non-lifting morphism}
 Let $Y=[-1,1]$ and consider the map
 $\rho:Y=[-1,1]\rightarrow \Bir(\mathbb P^2_\R)$, 
 \begin{equation*}
       \rho(t)\colon [x_0:x_1:x_2]\mapsrat
    \begin{cases}
    [x_0 \cdot P^t(x_0,x_1):x_1 \cdot P^t(x_0,x_1)
 :x_2\cdot P^t(x_0,x_1)+tx_0x_1], & t\neq 0 \\
      [x_0:x_1:x_2], &  t=0
    \end{cases}
  \end{equation*}
$$\text{ where  }\ P^t(x_0,x_1)=\cos\Big(\frac{2\pi}{t}\Big)x_0+ x_1.$$
Remark that locally, in the chart $x_0=1$, the map $\rho(t)$ is given by
$$(x_1,x_2)\mapsto \left(x_1,x_2+\frac{tx_1}{\cos(\frac{2\pi}{t})+x_1}\right).$$
The inverse of $p(t)$ is given by 
\[
\rho(t)^{-1}\colon [x_0:x_1:x_2] \mapsrat [x_0\cdot P^t(x_0,x_1):x_1\cdot P^t(x_0,x_1):x_2\cdot P^t(x_0,x_1)-tx_0x_1]
\]
i.e. $\rho(t)^{-1}=\rho(-t)$. 

We have $\rho(t)\mapsto \id$ if $t\mapsto 0$, $t\neq 0$. 
Indeed, $\rho(t)$ is of degree $2$ for any $t\neq0$, and we consider the set $A:=\pi_2^{-1}(\rho([-1,1]\setminus\{0\}))\subset H_2(\R)$. 
Since the (Euclidean) closure of the graph $\Gamma:=\{(t,\cos(\frac{2\pi}{t}))\mid t\neq0\}$ is equal to $\Gamma\cup\{0\}\times[-1,1]$, the (Euclidean) closure of $A$ in $H_2(\R)$ is
\[
\bar{A}=A\cup\{[x_0(sx_0+x_1):x_1(sx_0+x_1):x_2(sx_0+x_1)]\in H_2(\R)\mid s\in[-1,1]\}.
\]
We have $\pi_2(\bar{A}\setminus A)=\{\id\}$, so it follows that $\rho(t)\rightarrow\id$ for $t\mapsto 0$. 
Therefore, the map $\rho:[-1,1]\rightarrow \Bir(\mathbb P^2_\R)$ is continuous.

The maps $\rho(t)$, $t\neq 0$, are regular everywhere except at $[0:0:1]$
and they contract the line $\Big\{P^t(x_0,x_1)=\cos\Big(\frac{2\pi}{t}\Big)x_0+ x_1=0\Big\}$ to the point $[0:0:1]$.
Moreover, each $\rho(t)$ fixes $p=[1:0:0]$, i.e.\ $\rho(t)(p)=p$.
\end{example}

The following statement implies Proposition~\ref{prop:nolifting-main}: for $n\geq3$, Example~\ref{ex: non-lifting morphism} can be naturally extended to the continuous map 
$$[-1,1]\to\Bir(\p^n_\k),\quad t\mapsto(x_1,x_2+\frac{tx_1}{\cos(\frac{2\pi}{t})+x_1},x_3,\dots,x_n)$$
 and the same arguments apply.

\begin{proposition}
\label{prop:nolift}
Let $\k=\R$ or $\k=\C$.
Let $\rho\colon [-1,1]\to\Bir(\p^2_\k)$ be as in Example~\ref{ex: non-lifting morphism}.
For arbitrarily small $\varepsilon>0$ and any $d\geq1$, there does not exist a continuous map $\rho_{\varepsilon}\colon(-\epsilon,\epsilon)\to H_d(\k)$ such that $\rho|_{(-\epsilon,\epsilon)}=\pi_d\circ\rho_{\varepsilon}$.
\end{proposition}
\begin{proof}
Assume that there exists $d\geq1$ and $\varepsilon>0$ and a lift $\hat\rho:(-\epsilon,\epsilon)\rightarrow H_d(\k)$. 
Remember that $H_d(\k)\subset W_d(\k)$ where $W_d(\k)$ denotes the projective space parametrising equivalence classes of non-zero  $3$-tuples $(h_0,h_1,h_2)$
of homogeneous polynomials $h_i\in \k[x_0,x_1,x_2]$ of degree $d$,
where $(h_0,h_1,h_2)$ is equivalent to $(\lambda h_0,\lambda h_1,\lambda h_2)$ for any $\lambda\in \k^{*}$,
i.e.\ $W_d=\mathbb P(\k[x_0,x_1,x_2]_d\times \k[x_0,x_1,x_2]_d \times \k[x_0,x_1,x_2]_d)$ where
$\k[x_0,x_1,x_2]_d$ is the space of homogeneous polynomials of degree $d$.
The quotient map $\pi_{W_d}:(\k[x_0,x_1,x_2]_d\times \k[x_0,x_1,x_2]_d \times \k[x_0,x_1,x_2]_d)\setminus\{0\}\rightarrow W_d$ is a surjective submersion and 
thus locally there exist smooth sections of $\pi_{W_d}$. 
Therefore, (after possibly further shrinking $\Omega=(-\varepsilon, \varepsilon)$) the map $\rho$ can also be lifted to 
a continuous map $\tilde\rho:\Omega\rightarrow \k[x_0,x_1,x_2]_d\times \k[x_0,x_1,x_2]_d \times \k[x_0,x_1,x_2]_d$ with 
$\rho=\pi_d\circ\pi_{W_d}\circ \tilde\rho$. In other words, there are homogeneous polynomials $h_0^t$, $h_1^t$ and $h_2^t$ of degree $d$ whose coefficients depend 
continuously on $t$ and which satisfy $\rho(t)=\pi_d((h_0^t:h_1^t:h_2^t))$ for all $t\in \Omega =(-\varepsilon,\varepsilon)$.

In particular, for any $x_0\neq0$, $x_2$ and any $t\in (-\varepsilon,\varepsilon)\setminus\{0\}$ with $\frac{2\pi}{t}\notin \mathbb Z$, 
i.e.\ $\cos\left(\frac{2\pi}{t}\right)\neq 0$, we have 
$$\rho(t)([x_0:-\cos\left(\frac{2\pi}{t}\right)x_0 :x_2])= [0:0: -t\cos\left(\frac{2\pi}{t}\right)x_0^2] = [0:0:1]$$
and consequently,
$$h_0^t(x_0,-\cos\left(\frac{2\pi}{t}\right)x_0 ,x_2) = h_1^t(x_0,-\cos\left(\frac{2\pi}{t}\right)x_0 ,x_2) = 0$$
for all those $t$, $x_0$ and $x_2$.
For any $s\in[-1,1]$ there is a sequence $t_m \in (-\varepsilon,\varepsilon)\setminus\{0\}$ with $\frac{2\pi}{t_m}\notin \mathbb Z$ for all $m$
and $t_m \to 0$ such that $\cos\left(\frac{2\pi}{t_m}\right) \to s$.
Since the coefficients of the homogeneous polynomials $h_i^t$ continuously depend on $t$ we get
$$h_i^{0}(x_0,-sx_0,x_2) = \lim_{m\rightarrow\infty} h_i^{t_m}(x_0,-\cos\left(\frac{2\pi}{t_m}\right)x_0 ,x_2) = 0$$
for $i = 0,1$ and any $s\in [-1,1]$, $x_0\neq 0$ and $x_2$. 
Since the $h_i^t$ are all homogeneous polynomials, this implies $h_0^0 = h_1^0 = 0$ and hence $h_2^0$ has to be non-trivial and there exists a $q\in \k^3$ with 
$[q]\neq [0:0:1]$ in $\p^2_\k$ and $h_2^0(q) \neq 0$. 
Consequently, we get 
$\rho(0)([q]) =(h_0^0: h_1^0: h_2^0)[q] = [0:0:1],$
which is contradicting the original definition of $\rho(0) = \id$.
\end{proof}

\end{document}